\documentclass{article}
\usepackage[utf8]{inputenc}
\usepackage{hyperref}
\usepackage{microtype}
\usepackage{mathtools}
\usepackage{amssymb}
\usepackage{amsthm}

\usepackage{tikz}

\DeclarePairedDelimiter{\rawbrackets}{[}{]}
\newcommand{\brackets}{\rawbrackets*}
\DeclarePairedDelimiter{\rawabs}{\lvert}{\rvert}
\newcommand{\abs}{\rawabs*}
\DeclarePairedDelimiter{\rawp}{(}{)}
\newcommand{\p}{\rawp*}
\DeclarePairedDelimiter{\rawbraces}{\{}{\}}
\newcommand{\braces}{\rawbraces*}
\DeclarePairedDelimiter{\rawclopen}{[}{)}
\newcommand{\clopen}{\rawclopen*}
\newcommand{\LCS}{\operatorname{LCS}}
\newcommand{\SHIFT}{\operatorname{SHIFT}}
\newcommand{\BigSHIFT}{\operatorname{BigSHIFT}}
\newcommand{\eqdef}{\coloneqq}
\DeclareMathOperator{\Span}{span}
\newcommand{\E}{\mathbf{E}}
\renewcommand{\Pr}{\mathbf{Pr}}
\newcommand{\definition}{\textit}

\newtheorem{lem}{Lemma}
\newtheorem{thm}{Theorem}

\title{Length of the longest common subsequence between overlapping
words}
\author{Boris Bukh\thanks{Department of Mathematical Sciences, Carnegie
Mellon University, Pittsburgh, PA 15213, USA. Supported in part by Sloan
Research Fellowship and by U.S.\ taxpayers through NSF CAREER grant
DMS-1555149.} \and Raymond Hogenson\thanks{This research was
supported by the CMU SURF program and by NSF CAREER grant DMS-1555149.}}

\begin{document}
\maketitle
\begin{abstract}
Given two random finite sequences from $\brackets{k}^n$ such 
that a prefix of the first sequence is a suffix of the second,
we examine the length of their longest common subsequence. If
$\ell$
is the length of the overlap, we prove that the expected length of an
LCS is approximately $\max\p{\ell, \E\brackets{L_n}}$, where $L_n$ is the length of
an LCS between two independent random sequences. We also obtain tail
bounds on this quantity.
\end{abstract}

\section{Introduction}
A \definition{word} is a finite sequence of symbols over some alphabet.
We write $\abs{W}$ for
the length of a word~$W$. We write $W[i]$ for the $i$th symbol of $W$,
indexing starting with $0$.
A \definition{subsequence} of a word $W$ is a word obtained by deleting
symbols from $W$.
A \definition{common subsequence} between two words $V$ and $W$ is a
subsequence of both $V$ and $W$. A natural notion of similarity between
two words is the length of the \definition{longest common subsequence} (LCS) for the
two. We write $\LCS\p{V, W}$ for
the length of an LCS between words $V$ and $W$. A \definition{subword} of a
word $W$ is a subsequence consisting of contiguous symbols from $W$. We
denote the subword of $W$ consisting of symbols $a$ through $b - 1$ by
$W[a, b)$. For a set $A=\braces{0 \le i_1 < i_2 <
\dotsb}$, we
write $W[A]$ for the subsequence given by $W[i_1]W[i_2]\dotsm$.

To take a concrete example, let $V = \textsf{abedbba}$ and $W =
\textsf{aabdca}$. Then $\LCS\p{V, W} = 4$, as evidenced by the common
subsequence \textsf{abda}. In this paper we are interested 
in the LCS between words chosen randomly.

As the nature of the symbols will not be important to us, will use $\brackets{k}\eqdef \{1,2,\dotsc,k\}$ for the alphabet. 
We'll write $W \sim \brackets{k}^n$ to indicate that $W$ is a word chosen uniformly
at random from $\brackets{k}^n$.

Let $L_n \eqdef \LCS\p{V, W}$ where $V, W \sim \brackets{k}^n$. Then
we define
\begin{equation*}
\gamma_k
\eqdef \lim_{n \to \infty} \frac{\E\brackets{L_n}}{n}.
\end{equation*}
See~\cite{longest-common-subsequences-of-two-random-sequences} for a
proof that this limit exists, as well as
upper and lower bounds on $\gamma_k$.
Refer to~\cite{improved-bounds-on-the-average-length-of-longest-common-subsequences}
for the best known bounds on $\gamma_2$ and a deterministic method
to determine accurate bounds for $\gamma_k$ for $k > 2$.
In~\cite{expected-length-of-the-longest-common-subsequence-for-large-alphabets}
it is shown that $\gamma_k \sqrt{k} \to 2$ as $k \to \infty$.

In this paper, we examine a related problem: the LCS between two random
words which overlap. Namely,
let $0 \le \alpha \le 1$, pick $Z \sim \brackets{k}^{n
+ \alpha n}$, and choose $V = Z[0, n)$ and $W = Z[\alpha n, n + \alpha n)$. Thus a
suffix of $V$ is the same as a prefix of $W$. We say that \definition{$W$ is shifted from $V$ by $\alpha$}.
We will examine $\SHIFT\p{n, k, \alpha n} \eqdef \LCS\rawp[\big]{Z[0, n),
Z[\alpha n, n + \alpha n)}$ where $Z \sim \brackets{k}^{n + \alpha n}$.

This is motivated in part by an application to DNA sequencing. In this
process, we have two sections of DNA which can be regarded as words over
the alphabet of nucleotides. The
pieces of DNA may overlap, and we wish to determine whether the
similarity between them is more than coincidence, i.e., if they are
indeed from the same section of the genome.

\textbf{Acknowledgment:} We thank Mikl\'os R\'acz for telling us of the problem.

\section{Results}
For $W$ shifted from $V$ by $\alpha n$, the length of the
overlap ($\alpha n$) is a lower bound on the length of the LCS since the
overlapping section is \emph{a} common subsequence between the two. When
$\alpha n$ is much less than $\E[L_n]$,
we might think that the overlap does not matter and
$\LCS\p{W, V}$ behaves like $L_n$. This is indeed so, as the following
two theorems show.
\begin{thm}\label{upper}
There exists a constant $c_k$ such that for any $t
\ge 6\sqrt{n}$,
\begin{equation*}
\Pr\brackets{\SHIFT\p{n, k, \alpha n} \ge \max\p{n - \alpha n + 1, \gamma_kn +
t}} \le \exp\p{-c_kt^2/n}
\end{equation*}
when $n$ is sufficiently large.
\end{thm}

\begin{thm}\label{lower}
There exists a constant $c_k$ such that for any $t \ge 5n^{3/4}\sqrt{\log n}$,
\begin{equation*}
\Pr\brackets{\SHIFT\p{n, k, \alpha n} \le \gamma_k n - t} \le
\exp\rawp[\big]{-c_kt^2/n^{3/4}}
\end{equation*}
when $n$ is sufficiently large.
\end{thm}

All logarithms in this paper are to base $e = 2.71\dots$.

We expect that $\SHIFT\p{n, k, \alpha n} = \max\rawp[\big]{n - \alpha n, \E[L_n] +
O\p{\sqrt{n}}}$ with high probability. This is supported
by~\cite{standard-deviation-of-the-longest-common-subsequence} which
shows that the standard deviation of $L_n$ is $O\p{\sqrt{n}}$.

\section{Tools}
Here we collect several auxiliary results.

\begin{lem}[\cite{the-rate-of-convergence-of-the-mean-length-of-the-longest-common-subsequence},
Theorem~1.1]\label{rateofconvergence}
\begin{equation*}
\gamma_k n \ge \E\brackets{L_n} \ge \gamma_k n - 4\sqrt{n \log n}.
\end{equation*}
\end{lem}

Let $\Omega = \prod_{i = 1}^n \Omega_i$ where each $\Omega_i$ is a
probability space and $\Omega$ has the product measure. Let $h \colon
\Omega \to \mathbb{R}$. Let $X$ be a random variable given by $X =
h\p{\cdot}$.

We call $h \colon \Omega \to \mathbb{R}$ \definition{Lipschitz} if
$\abs{h\p{x} - h\p{y}} \le 1$ whenever $x, y$ differ in at most one
coordinate.

\begin{lem}[Azuma's inequality, \cite{the-probabilistic-method},
Theorem~7.4.2]\label{azuma}
If $h$ is Lipschitz, then
\begin{equation*}
\Pr\brackets{\abs{X - \E[X]} > \lambda } \le e^{-t^2/4}.
\end{equation*}
\end{lem}

\begin{lem}[Hoeffding's inequality,
\cite{probability-inequalities-for-sums-of-bounded-random-variables}]\label{hoeffding}
Suppose $X_i$ are independent random variables with $a_i \le X_i \le
b_i$. Then for all $t > 0$,
\[\Pr\brackets{\sum_i X_i \le \E\brackets{X} - t} \le \exp\p{\frac{-2t^2}{\sum_i
\p{b_i - a_i}^2}}.\]
\end{lem}

\begin{lem}\label{boris}
Let $X_1,\dotsc,X_m$ be independent random variables, each of which is
exponential with mean~$k$. Let $X=X_1+\dotsb+X_m$. Then
\begin{equation*}
\Pr\brackets{X\leq m\p{k-\lambda}} \le \exp\p{-m\lambda^2/2k^2}.
\end{equation*}
\end{lem}
\begin{proof}
We may assume that $\lambda<k$, for otherwise the result is trivial.
With hindsight set $t=\lambda/k\p{k-\lambda}$.
Then
\begin{align*}
  \E\brackets{\exp(-tX_i)}&=\sum_{j\geq 1} \exp\p{-tj}
  \p{1-1/k}^{j-1}\p{1/k}\\
               &=\frac{1}{1+\p{e^t-1}k}\leq \frac{1}{1+kt}.
\end{align*}
From this it follows that $\E\brackets{\exp\p{tmk-tX}} \le \p{e^{tk}/\p{1+kt}}^m$.
Therefore by Markov's inequality we have
\begin{align*}
  \Pr\brackets{mk-X\geq \lambda m}&\leq
  \frac{\p{e^{tk}/(1+kt)}^m}{\exp\p{t\lambda
  m}}=\p{\frac{e^{t(k-\lambda)}}{1+kt}}^m\\
&=\p{e^{\lambda/k}(1-\lambda/k)}^m\leq
\p{1-\p{\lambda/k}^2/2}^m\\
&\leq \exp\p{-m\lambda^2/2k^2}
\end{align*}
where in the penultimate line we used the inequality $e^x\p{1 - x} \le 1
- x^2/2$, which can be established by considering the Taylor expansion
of $e^x\p{1 - x}$.
\end{proof}

In this paper, we will weaken this bound to
\begin{equation}\label{boriseasy}
\Pr\brackets{X \le m\p{k - \lambda}} \le \exp\p{-m\lambda/4k + m/32}
\end{equation}
using $x^2 \ge x/2 - 1/16$.

\section{Proof of Theorem~\ref{upper}}
\begin{figure}
\centering
\begin{tikzpicture}
\draw (0, 0) node [left] {$V$} -- (5, 0);
\draw (3, 1) node [left] {$W$} -- (8, 1);

\def\oa{(1,0)}
\def\ob{(3.5,0)}
\def\oc{(3.75,0)}
\def\od{(4,0)}
\def\oe{(4.5,0)}

\def\ta{(3.25,1)}
\def\tb{(4,1)}
\def\tc{(4.5,1)}
\def\td{(6,1)}
\def\te{(8,1)}

\draw \oa -- \ta;
\draw \ob -- \tb;
\draw \oc -- \tc;
\draw \od -- \td;
\draw \oe -- \te;

\filldraw \oa circle [radius=2pt]
                 \ob circle [radius=2pt]
                 \oc circle [radius=2pt]
		 \od circle [radius=2pt]
		 \oe circle [radius=2pt]
		 \ta circle [radius=2pt]
		 \tb circle [radius=2pt]
		 \tc circle [radius=2pt]
		 \td circle [radius=2pt]
		 \te circle [radius=2pt];
\end{tikzpicture}
\caption{A visual representation of a common subsequence between shifted
words.\label{mental-model}}
\end{figure}
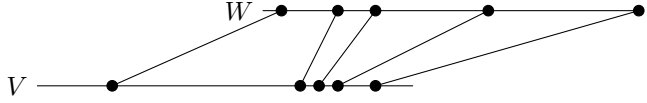

\begin{figure}
\centering
\begin{tikzpicture}
\draw (0, 0) -- (5, 0);
\draw (2, 1) -- (7, 1);

\def\oa{(2,0)}
\def\ob{(2.5,0)}
\def\oc{(3,0)}
\def\od{(3.5,0)}
\def\oe{(4,0)}
\def\of{(4.5,0)}
\def\og{(5,0)}

\def\ta{(2,1)}
\def\tb{(2.5,1)}
\def\tc{(3,1)}
\def\td{(3.5,1)}
\def\te{(4,1)}
\def\tf{(4.5,1)}
\def\tg{(5,1)}

\draw \oa -- \ta;
\draw \ob -- \tb;
\draw \oc -- \tc;
\draw \od -- \td;
\draw \oe -- \te;
\draw \of -- \tf;
\draw \og -- \tg;

\filldraw \oa circle [radius=2pt]
                 \ob circle [radius=2pt]
                 \oc circle [radius=2pt]
		 \od circle [radius=2pt]
		 \oe circle [radius=2pt]
		 \of circle [radius=2pt]
		 \og circle [radius=2pt]
		 \ta circle [radius=2pt]
		 \tb circle [radius=2pt]
		 \tc circle [radius=2pt]
		 \td circle [radius=2pt]
		 \te circle [radius=2pt]
		 \tf circle [radius=2pt]
		 \tg circle [radius=2pt];
\end{tikzpicture}
\caption{An example where the common subsequence is the overlapping
section.\label{mental-overlap}}
\end{figure}
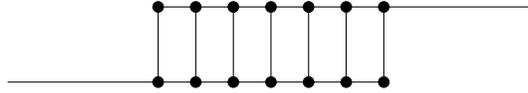

There is a geometric way to interpret a common subsequence. Consider a line segment
from $\p{0, 0}$ to $\p{n - 1, 0}$, and a second from $\p{\alpha n, 1}$ to $\p{n +
\alpha n - 1, 1}$. Now we place the symbols from $V$ on the first line segment
and those from $W$ on the second. For each pair of symbols in $V$ and
$W$, connect them with an edge if the symbols are equal. The LCS, then,
is the largest set of noncrossing edges. Furthermore, symbols
aligned vertically will be certainly equal by the nature of the shift.
See Figs.~\ref{mental-model}
and~\ref{mental-overlap}
for examples.

\begin{lem}\label{grosscase}
Suppose we have indices $0 \le i_1 < \dots < i_t < n$.
Let $Z \sim \brackets{k}^{n + \alpha n}$ and take $V = Z[0,n)$ and $W =
Z\clopen{\alpha n, n + \alpha n}$.
Define $A =
V\brackets{\braces{i_1, \dots, i_t}}$, $B_{\text{start}} = \max\p{i_t
- t + 1 - \alpha n, 0}$, and $B = W\clopen{B_{\text{start}}, n}$. Then $A\brackets{\ell}$ is to the left of
$B\brackets{\ell}$ for every $\ell$, and
\begin{equation*}
\Pr\brackets{\LCS\p{A, B} = t} \le \exp\rawp[\big]{\p{\abs{B} -
\tfrac{7}{8}k\abs{A}}/4k}
\end{equation*}
\end{lem}
\begin{proof}
See Fig.~\ref{lemmasubword}
\begin{figure}
\centering
\begin{tikzpicture}
\draw (0,0) node [left] {V} -- (7,0);
\draw (3,1) node [left] {W} -- (10,1);
\filldraw (5.5,0) node {} circle [radius=2pt];
\filldraw (6,0) node {} circle [radius=2pt];
\filldraw (6.2,0) node [below=0.8ex] {A} circle [radius=2pt];
\filldraw (6.5,0) node {} circle [radius=2pt];
\filldraw (6.8,0) node {} circle [radius=2pt];
\draw [very thick] (6,1) -- node [above=0.5ex] {B} (10,1);
\end{tikzpicture}
\caption{The dots represent the symbols of $V$ that are in $A$ and the bold
line represents the symbols in $B$.\label{lemmasubword}}
\end{figure}
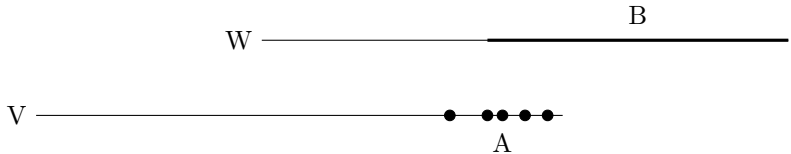
for a depiction of this lemma in our geometric model.
We will use the notation $W_1 \trianglelefteq W_2$ to mean that $W_1$ is a
subsequence of $W_2$.

We prove first that $A\brackets{\ell}$ is
always strictly to the left of $B\brackets{\ell}$ for every $\ell$. Let
$\ell_V = i_\ell$, and $\ell_W = \ell + B_{\text{start}} \ge \ell + i_t
- t + 1 - \alpha n$. These
are the positions of $A[\ell]$ and $B[\ell]$ in $V$ and $W$ respectively.
The horizontal position of $\ell_W$ is then at least
$\ell + i_t - t + 1$ (from the left of $V$), and the position of $\ell_V$
is at most $i_t - t + \ell$.\medskip

To prove the bound on $\Pr\brackets{\LCS\p{A,B}=t}$, we introduce an equivalent
way of generating random word $Z$. Let $R, S \sim \brackets{k}^\infty$.
Imagine $n + \alpha n$ placeholders corresponding to symbols of $Z$. In
the beginning, the placeholders are empty. We will use symbols from $S$
and from $R$ in order to fill the placeholders using the following
process. Start with $\ell_A = \ell_B = 0$. At each step, if the
placeholder for $A\brackets{\ell_A}$ is empty, use the next symbol from
$S$ to fill it. Then we examine successive symbols from $R$ until the
last examined symbol is equal to $A\brackets{\ell_A}$; we use the
examined symbols to fill placeholders in $B$ starting from
$B\brackets{\ell_B}$ (if we run out of empty placeholders, we simply
discard symbols from $R$). Finally, we increment $\ell_A$ and increase
$\ell_B$ appropriately so that $B\brackets{\ell_B}$ is the first
unfilled placeholder in $B$. 

Note that at each step $\ell_B$ increases by at least $1$. Since $\ell_A$ 
increases by exactly $1$, it follows that $B\brackets{\ell_B}$ is to the right of
$A\brackets{\ell_A}$ at all times in this process.

Finally, after all placeholders in $A$ are filled, we fill the rest of $Z$ with
symbols from $S$. 

During this process, each next filled symbol is independent of all
the ones before. Therefore, the word $Z$ we obtain is a uniformly random word.

Let $X_i$ be the number of symbols from $R$
consumed after we match the $\p{i - 1}$th symbol of $A$ but before we match
the $i$th symbol of $A$. Let $X = \sum_i X_i$.
Then $\Pr\brackets{A
\trianglelefteq B} = \Pr\brackets{X \le \abs{B}}$.

We apply~\eqref{boriseasy} (which is a weakening of Lemma~\ref{boris}) with $m = \abs{A}$, $\lambda = k -
\abs{B}/\abs{A}$:
\begin{equation*}
\Pr\rawbrackets[\big]{A \trianglelefteq B} =
\Pr\brackets{X \le \abs{B}} \le
\exp\p{\p{\abs{B} - \tfrac{7}{8}k\abs{A}}/4k}.
\qedhere
\end{equation*}
\end{proof}

Pick
$Z \sim \brackets{k}^{n + \alpha n}$, and $V = Z[0, n)$, $W = Z[\alpha n, n +
\alpha n)$ (so $W$ is shifted from $V$ by $\alpha n$).

We define the \definition{span} of an edge $e$ between the $i$th symbol from $V$ and the $j$th
symbol from $W$ to be $\Span\p{e}\eqdef j + \alpha n - i$. This is the
difference in $x$-coordinates in the geometric model.
If the span is positive, the slope of the edge
will be positive, and conversely a negative span indicates a
negative-sloping edge. In particular, if the span is $0$, the symbols
will be equal due to the nature of the shift. We say that symbols are
overlapping if the span between them is $0$.

\medbreak

With hindsight, set $\varepsilon = 0.01$. We will break our analysis
into several cases, according to the shape of the LCS. For each shape,
we'll bound the probability that there is a long LCS between $V$ and $W$
of that shape. Each case can be described by the edge with the least
span. Let $e$ be an edge connecting a symbol of~$V$ with a symbol of~$W$,
and define the following random variable
\begin{center}
$\BigSHIFT_e\eqdef\null$``Maximum length of a common subsequence between $V$ and $W$ 
that uses the edge $e$, and the edge $e$ is an edge of largest span in this subsequence''. 
\end{center}

Note that the probability of the event in Theorem~\ref{upper} can be bounded by
\begin{equation}\label{unionbound}
\begin{aligned}
  &\Pr\brackets{\SHIFT\p{n, k, \alpha n} \ge \max\p{n - \alpha n + 1, \gamma_kn +
t}}\\&\qquad\le \sum_{e} \Pr\brackets{\BigSHIFT_e\ge \max\p{n-\alpha n+1,\gamma_k n + t}},
\end{aligned}
\end{equation}
where the sum is over all edges~$e$ connecting a symbol of $V$ with a symbol of~$W$.

Let $e$ be an arbitrary edge, connecting $V\brackets{i}$ and $W\brackets{j}$. We shall estimate
the $\Pr\brackets{\BigSHIFT_e\ge \max\p{n-\alpha n+1,\gamma_k n + t}}$.

\paragraph{Case $\Span\p{e} \le 0$.} There are $n
- i$ symbols to the right of $e$ in $V$ and $j - 1$ symbols to the left
of $e$ in $W$. Therefore the length of the LCS is at most $n - i +
j \le n - \alpha n$ so the event $\BigSHIFT_e$ cannot
occur.

\paragraph{Case $0 < \Span\p{e} \le \varepsilon n$.}
Write
$A = V[i + 1,n)$, $B = W[0,j)$, so that $A$ is the subword of $V$
after the symbol in position $i$ and $B$ is the subword of $W$
before position $j$.
Let $s = j + \alpha n - i$.
We have $s \le \varepsilon n$.

If $W_2$ is a subsequence of $W_1$,
then there are indices $i_1, \dots, i_{\abs{W_2}}$ such that $W_2[\eta]
= W_1[i_\eta]$ for each $1 \le
\eta \le \abs{W_2}$. We use the notation
$W_1 \setminus W_2$ to mean $W_1$ with symbols $i_1, \dots,
i_{\abs{W_2}}$ deleted.

Because $\abs{A} + \abs{B} = n - \alpha n + s$, for the length of the LCS to be greater
than or equal to $n - \alpha n$, we must have
\begin{equation*}
\text{$A \setminus L_1 \trianglelefteq W[j + 1, n)$ and $B \setminus L_2
\trianglelefteq V[0, i)$}
\end{equation*}
for some $L_1 \trianglelefteq A$, $L_2 \trianglelefteq B$ satisfying
$\abs{L_1} + \abs{L_2} = s$.

We can then bound $\Pr\brackets{\SHIFT\p{n, \alpha, k} > n - \alpha n}$ from
above by
\begin{equation}\label{almostverticalprbound}
\sum_{\substack{L_1 \subseteq A \\ L_2 \subseteq B \\ \abs{L_1} +
\abs{L_2} = s}}
\Pr\rawbrackets[\big]{\text{$A \setminus L_1 \trianglelefteq W[j + 1,n)$ and $B
\setminus L_2 \trianglelefteq V[0,i)$}}.
\end{equation}

\begin{figure}
\centering
\begin{tikzpicture}
\draw (0,0) node [left] {$V$} -- (6,0);
\draw (3,1) node [left] {$W$} -- (9,1);

\draw (4,0) -- node [left] {$e$} (4.5,1) node [above] {$W\brackets{j}$};

\draw [very thick] (4,0) node [below] {$V\brackets{i}$} -- node [above=0.4ex] {$A$} (6,0);
\draw [very thick] (3,1) -- node [above=0.4ex] {$B$} (4.5,1);
\end{tikzpicture}
\caption{Illustration of the case when $0 < \Span\p{e} \le \varepsilon
n$. The words $V$ and $W$ are divided by the edge which is assumed to
exist in this case.\label{randomizedprocess}}
\end{figure}
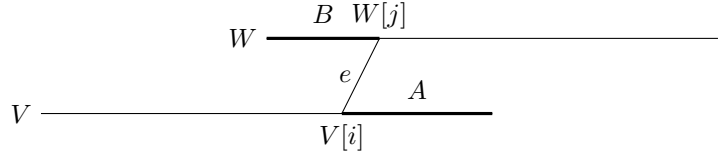

Let $A' = A \setminus L_1$ and $W' = W[j + 1, n)$.
Note that in our geometrical interpretation of the LCS, we have $A'$
positioned above $W'$, and slightly to the left. See
Fig.~\ref{randomizedprocess}.

Consider the starting position of $W'$. It is $j + 1$ into $W$. $j$ is
$s$ symbols right of $i$. So $W'$ is $i + s + 1$ from the beginning of
$V$. In the terminology of
Lemma~\ref{grosscase}, $i_t = i + \abs{A}$ and $t = \abs{A} - s + 1$, so
$i_t - t + 2 = i + s + 1$, and we can thus apply Lemma~\ref{grosscase}
to $A'$ and $W'$.
So we have
\begin{equation}\label{abound}
\Pr\brackets{A \setminus L_1 \trianglelefteq W[j + 1,n)} \le
\exp\p{\p{n - \abs{B} - \tfrac{7}{8}k\p{\abs{A} - \abs{L_1}}}/4k}.
\end{equation}
By the same reasoning
\begin{equation}\label{bbound}
\Pr\brackets{B \setminus L_2
\trianglelefteq V[0, i)} \le \exp\p{\p{n - \abs{A} - \tfrac{7}{8}k\p{\abs{B} -
\abs{L_2}}}/4k}.
\end{equation}
Therefore we combine~\eqref{almostverticalprbound}, \eqref{abound},
and~\eqref{bbound} to get
\begin{align}
\begin{split}
&\Pr\brackets{\BigSHIFT_e > n - \alpha n} \\
&\le \sum_{\substack{L_1 \subseteq A \\ L_2 \subseteq B \\ \abs{L_1} +
\abs{L_2} = s}}
\exp\p{\p{2n - \abs{A} - \abs{B} - \tfrac{7}{8}k\p{\abs{A} + \abs{B} -
\abs{L_1} - \abs{L_2}}}/4k}
\end{split} \\
&\le \sum_{\substack{L_1 \subseteq A \\ L_2 \subseteq B \\ \abs{L_1} + \abs{L_2} =
s}} \exp\p{\p{1 + \alpha - \varepsilon - \tfrac{7}{8}k + \tfrac{7}{8}\alpha k}n/4k}
\label{boundsimplify}
\end{align}
Where we simplified~\eqref{boundsimplify} using $\abs{A} + \abs{B} = n -
\alpha n + s$ and~$\abs{L_1} + \abs{L_2} = s$.

Now we simplify the last sum as
\begin{align}
\begin{split}\label{secondcasebound}
\eqref{boundsimplify} &= \binom{\abs{A} +
\abs{B}}{s}\exp\p{\p{1 + \alpha - \varepsilon - \tfrac{7}{8}k +
\tfrac{7}{8}\alpha
k}n/4k} \\
&\le \varepsilon n \exp\p{\varepsilon n \log\p{e\p{1 - \alpha +
\varepsilon}/\varepsilon} + \p{1 + \alpha - \varepsilon - \tfrac{7}{8}k +
\tfrac{7}{8}\alpha
k}n/4k} \\
&\le \exp\p{\p{1 + \alpha - \tfrac{7}{8}k + \tfrac{7}{8}\alpha k + \varepsilon - \varepsilon
\log \varepsilon}n/4k}
\end{split}
\end{align}
using $\binom{n}{k} \le \p{\frac{en}{k}}^k$.

\paragraph{Sub-case $n - \alpha n + \varepsilon n \le \gamma_k n$.}
An edge with span less than $\varepsilon n$ limits the
length of the LCS to
$n - \alpha n + \varepsilon n \le \gamma_k n$, so the event
$\BigSHIFT_e\ge \gamma_k n + t$ does not occur, i.e., the probability is $0$ in
this case.

\paragraph{Sub-case $n - \alpha n + \varepsilon n > \gamma_k n$.}
We can bound~\eqref{secondcasebound} above by
\begin{align*}
&\exp\p{\p{1 + \p{1 + \varepsilon - \gamma_k} - \tfrac{7}{8}k +
\tfrac{7}{8}\p{1 + \varepsilon -
\gamma_k}k + \varepsilon - \varepsilon \log \varepsilon}n/4k} \\
&\qquad \le \exp\p{\p{2 - \p{\tfrac{7}{8}k + 1 - \varepsilon}\gamma_k + 2\varepsilon -
\varepsilon \log \varepsilon}n/4k}
\end{align*}
since $\alpha \le 1 + \varepsilon - \gamma_k$ in this case.\medskip

Since there are no more than $n^2$ choices for edge $e$, the
contribution of this case to the right-hand side of \eqref{unionbound} 
is at most
\begin{align*}
&n^2\exp\p{\p{2 - \p{\tfrac{7}{8}k + 1 - \varepsilon} \gamma_k + 2\varepsilon -
\varepsilon \log \varepsilon}n/4k} \\
&\qquad \le \exp\p{\p{2 - \p{\tfrac{7}{8}k + 1 - \varepsilon} \gamma_k + 3 \varepsilon -
\varepsilon \log \varepsilon}n/4k}.
\end{align*}

We claim that $2 - \p{\frac{7}{8}k + 1 - \varepsilon}\gamma_k + 3\varepsilon -
\varepsilon \log \varepsilon < 0$. Recall that earlier we set $\varepsilon =
0.01$. So we must show that $2.08 - \p{\frac{7}{8}k + 0.99}\gamma_k <
0$. For $k \ge 3$ the bound $\gamma_k \ge 1/\sqrt{k}$
from~\cite{some-limit-results-for-longest-common-subsequences} suffices.
For $k = 2$, the lower bound on
$\gamma_k$
from~\cite{improved-bounds-on-the-average-length-of-longest-common-subsequences}
suffices to
show the inequality.
Therefore our upper bound in this case is
\begin{equation*}
\exp\p{-c_1n}
\end{equation*}
for some positive constant $c_1$ depending only on~$k$.

\paragraph{Case $\Span\p{e} > \varepsilon n$.}
To bound the probability of $\BigSHIFT_e$ in this case, we will
estimate the probability that there is a common subsequence 
of a certain approximate shape.

Formally, we say that a pair $(A_i,B_i)$ is a \definition{block}
if $A_i$ is a subword of $V_i$ and $B_i$ is a subword of $W_i$.
We say that $(A_i,B_i)_i$ is a \definition{block partition} 
if $(A_i)_i$ is a partition of $V$ and $(B_i)_i$ is a partition
of $W$. In this case, we call the edges between symbols of $A_i$ and $B_i$
for some $i$ \definition{dominated} by the partition. We say that a common
subsequence $C$ of $V$ and $W$ is dominated by the partition if all the edges
in geometric model are dominated by the partition.

A block partition $(A_i,B_i)_i$ is said to be \definition{nonoverlapping} if,
for each $i$, the symbols in $A_i$ and $B_i$ are disjoint. \medskip

\begin{figure}
\centering
\begin{tikzpicture}
\def\oa{(0,0)}
\def\ob{(1.5,0)}
\def\oc{(3.5,0)}
\def\od{(6.5,0)}
\def\oe{(7.5,0)}

\def\ta{(4,1)}
\def\tb{(5,1)}
\def\tc{(7,1)}
\def\td{(8,1)}
\def\te{(10,1)}

\draw[fill=lightgray,draw=none] (0,0) -- (2,0) -- \tc -- (3,1) -- cycle;
\draw[fill=lightgray,draw=none] (4,0) -- (6,0) -- \td -- cycle;

\draw (0, 0) -- (8, 0);
\draw (3, 1) -- (11, 1);

\draw[dashed] (0,0) -- (3,1);
\draw[dashed] (2,0) -- \tc;
\draw[dashed] (4,0) -- \td;
\draw[dashed] (6,0) -- \td;
\draw[dashed] (8,0) -- (11,1);

\draw (0, -0.1) -- (0, 0.1);
\draw (2, -0.1) -- (2, 0.1);
\draw (4, -0.1) -- (4, 0.1);
\draw (6, -0.1) -- (6, 0.1);
\draw (8, -0.1) -- (8, 0.1);

\draw \oa -- \ta;
\draw \ob -- \tb;
\draw \oc -- \tc;
\draw \od -- \td;
\draw \oe -- \te;

\filldraw \oa circle [radius=2pt]
                 \ob circle [radius=2pt]
                 \oc circle [radius=2pt]
		 \od circle [radius=2pt]
		 \oe circle [radius=2pt]
		 \ta circle [radius=2pt]
		 \tb circle [radius=2pt]
		 \tc circle [radius=2pt]
		 \td circle [radius=2pt]
		 \te circle [radius=2pt];
\end{tikzpicture}
\caption{\label{blocks}Edges from a common subsequence are shown by
solid lines and the block boundaries determined by these are shown by
dashed lines. Alternating blocks are also shaded for easy viewing.}
\end{figure}
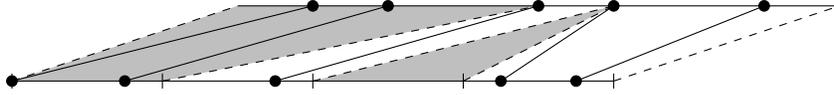

Given a common subsequence $C$ in which every edge satisfies $\Span\p{E}>\varepsilon n$,
we construct a partition that dominates it as follows.
We first partition $V$ into subwords of length exactly $\varepsilon n$. Let $(V_i)$ be this
partition. Consider the first edge in $S$ from $V_i$ to $W$. 
Let $W[S_i]$ be its endpoint in~$W$. Then $W$ is partitioned into subwords
$W_i\eqdef [S_i,S_{i+1})$. See Fig.~\ref{blocks}
for an example.

Because $\Span\p{e} > \varepsilon n$ for every edge $e$ and each $V_i$
has length $\varepsilon n$, we see that $V_i$ and $W_i$ do not overlap.
Hence, the resulting block partition is nonoverlapping.\medskip

So, it suffices to bound that probability that, given a nonoverlapping
partition $(V_i,W_i)_i$ into $1/\varepsilon$ blocks, there is a long common
subsequence dominated by this partition. The key observation
is that because $V_i$ and $W_i$ do not overlap the symbols
in $V_i$ and $W_i$ are independent.

Let $V'$ and $W'$ be two random words of length $n$ each, which are
disjoint from one another and from $Z$. Partition $V'$ into $(V_i')_i$
and $W'$ into $(W_i')_i$ in such a way that $\abs{V_i'}=\abs{V_i}$
and $\abs{W_i'}=\abs{W_i}$. Consider 
\begin{align*}
X&=\sum_i \LCS(V_i,W_i),\\
X'&=\sum_i \LCS(V_i',W_i').
\end{align*} 
Since $\E[\LCS(V_i,W_i)]=\E[\LCS(V_i',W_i')]$, by the linearity of expectation it follows 
that $\E[X]=\E[X']$. It is also clear that $X\leq \LCS(V',W')$, implying that
$\E[X']\leq \E[L_n]$. Therefore, $\E[X]\leq \E[L_n]$. 

From Azuma's inequality (Lemma~\ref{azuma}) and the bound on $\E[L_n]$ in Lemma~\ref{rateofconvergence} 
we obtain that
\[
  \Pr\brackets{X\geq \gamma_k n + t }\leq \Pr\brackets{X\geq \E[L_n] + t } \leq \exp\p{-t^2/4n}.
\]

Random variable $X$ is the longest length of a common subsequence that is dominated by a given $(V_i,W_i)_i$.
Since that there no more than $\binom{n}{1/\varepsilon}^2$ block partitions into $1/\varepsilon$ blocks,
from the union bound it follows that the contribution of the case $\Span\p{e} > \varepsilon n$ to \eqref{unionbound} is at most
\[
  \binom{n}{1/\varepsilon}^2\exp\p{\frac{-t^2}{4n}} \le \exp\p{\frac{-t^2}{5n}}
\]
since $\varepsilon = 0.01$ and $n$ is large enough.

Summing all the contributions from all the cases, we see that 
the right side of \eqref{unionbound} is bounded by
\begin{equation*}
\exp\p{-c_1n} + \exp\p{-c_2t^2/n} \le \exp\p{-c_kt^2/n}
\end{equation*}
for some constant $c_k$ as long as $n$ is large enough.
\qed

\section{Proof of Theorem~\ref{lower}}
We have two words $V$ and $W$ shifted by $\alpha n$. Divide each word
into blocks of size $\alpha n^{1/2}$. Write $V_i$ for the $i$th block, $V_i =
V\clopen{i\alpha n^{1/2}, \p{i + 1}\alpha n^{1/2}}$.
Similarly write $W_i$ for the $i$th block from $W$.

Note that $\LCS\p{V, W} \ge \sum_i \LCS\p{V_i, W_i}$,
therefore we can bound
\begin{equation}\label{lowerpreliminary}
\Pr\rawbrackets[\big]{\LCS\p{V, W} \le \gamma_k n - \ell} \le
\Pr\brackets{\sum_i \LCS\p{V_i, W_i} \le \gamma_k n - \ell}.
\end{equation}
Applying Lemma~\ref{rateofconvergence}, we have that
\begin{align*}
\gamma_k n &= \frac{\sqrt{n}}{\alpha}\gamma_k \alpha\sqrt{n} \\
&\le \frac{\sqrt{n}}{\alpha}\p{\E\brackets{L_{\alpha\sqrt{n}}} +
4\sqrt{\alpha\sqrt{n}\log\p{\alpha\sqrt{n}}}} \\
&\le \E\brackets{\frac{\sqrt{n}}{\alpha} L_{\alpha\sqrt{n}}} +
4n^{3/4}\sqrt{\log n}
\end{align*}

Thus we can upper-bound~\eqref{lowerpreliminary} by
\begin{equation}\label{readyforchernoff}
\Pr\brackets{\sum_i \LCS\p{V_i, W_i} \le
\E\brackets{\sum_i \LCS\p{V_i, W_i}} + 4n^{3/4}\sqrt{\log n} - t}.
\end{equation}
To bound~\eqref{readyforchernoff} we apply
Lemma~\ref{hoeffding}:
\begin{equation*}
\exp\p{\frac{-2\p{t - 4 n^{3/4}\sqrt{\log
n}}^2}{n^{1/2}\p{\alpha n^{1/2}}^2/\alpha }} \le \exp\p{\frac{-c_kt^2}{n^{3/4}}}
\end{equation*}
for an appropriately small $c_k$ if $n$ is large.
\qed

\bibliographystyle{abbrv}
\bibliography{bibliography}
\end{document}